\documentclass[10pt]{article}
\usepackage{amsfonts}
\usepackage{enumerate}
\usepackage{amsmath,amsthm, amssymb}
\usepackage[all,arc]{xy}
\title{Hausdorff Dimension of Average Conformal Hyperbolic Sets}
\author{Paul Wright}
\date{\today}

\newtheorem{thm}{Theorem}[section]
\newtheorem{cor}[thm]{Corollary}
\newtheorem{prop}[thm]{Proposition}
\newtheorem{lem}[thm]{Lemma}
\newtheorem{defn}[thm]{Definition}

\newtheorem{rem}[thm]{Remark}

\newcommand{\bpm}{\begin{pmatrix}}
\newcommand{\epm}{\end{pmatrix}}

\begin{document}
\maketitle

\begin{abstract}
The Hausdorff dimension of a conformal repeller or conformal hyperbolic set is well understood. For non-conformal maps, the Hausdorff dimension is only known in some special cases. Ban, Cao and Hu defined the concept of an average conformal repeller which generalises conformal, quasi-conformal and weakly conformal repellers, and they found an equation for the Hausdorff dimension for an average conformal repeller. In this paper we generalise this concept to average conformal hyperbolic sets, and obtain a similar equation for the Hausdorff dimension.
\end{abstract}

\section{Introduction}

The dimension of invariant sets such as repellers and hyperbolic sets is an important topic in dynamical systems. In the case of conformal repellers and conformal hyperbolic sets, the Hausdorff dimension and lower and upper box dimensions all agree, and are equal to the root of the so-called Bowen's equation. It is still an open problem to find the dimension of the invariant set of a non-conformal map, although some progress has been made (see \cite{survey}). In \cite{average conformal 1}, the authors introduced the concept of an \textit{average conformal repeller}, which for $C^1$ maps is more general than conformal, quasi-conformal and weakly conformal repellers. They obtained an equation for the Hausdorff dimension of an average conformal repeller. In this paper, we generalise quasi-conformal and average conformal maps to hyperbolic sets. We sketch a proof that for $C^1$ maps, average conformal hyperbolic sets are more general than quasi-conformal hyperbolic sets. Then we obtain an equation for the dimension of an average conformal hyperbolic set using similar arguments to \cite{average conformal 1}. Some improvements have been made to the structure and details of the proofs. 

Let $M$ be a compact Riemannian manifold, and let $f: M \rightarrow M$ be a $C^1$ diffeomorphism with a compact invariant set $\Lambda$. Then $\Lambda$ is a repeller if all the Lyapunov exponents of $f|_{\Lambda}$ are positive. The invariant set $\Lambda$ is called a \textit{hyperbolic set} if there exists a continuous splitting of the tangent bundle $T M = E^{(s)} \oplus E^{(u)}$, and constants $C > 0$, $0 < \lambda < 1$ such that for every $x \in \Lambda$, 
\begin{enumerate}
	\item $d_xf E^{(s)}(x) = E^{(s)}(f(x))$, $d_xf E^{(u)}(x) = E^{(u)}(f(x))$.
	\item For all $n \geq 0$, $\|df^n v\| \leq C \lambda^n \|v\|$ if $v \in E^{(s)}(x)$, and $\|df^{-n} v\| \leq C \lambda^n \|v\|$ if $v \in E^{(u)}(x)$.
\end{enumerate}
At each point $x \in \Lambda$ there are local stable and unstable manifolds $W^{(s)}(x)$ and $W^{(u)}(x)$. A hyperbolic set $\Lambda$ is called \textit{locally maximal} if there exists a neighbourhood $U$ of $\Lambda$ such that for any closed $f$-invariant subset $\Lambda'$ we have $\Lambda' \subseteq U$.

In \cite{average conformal 1}, a repeller is called \textit{average conformal} if there is exactly one unique Lyapunov exponent with respect to any invariant measure. In their formulation, there are several Lyapunov exponents, all of which are equal; in this paper, we say there is one Lyapunov exponent. In this paper, a hyperbolic set $\Lambda$ will be called \textit{average conformal} if it has two unique Lyapunov exponents, one positive and one negative. That is, for any invariant measure $\mu$, the Lyapunov exponents are $\chi^{(s)}(\mu) < 0 < \chi^{(u)}(\mu)$. For convenience we will also say that the map $f$ is average conformal if it has an average conformal repeller or an average conformal hyperbolic set.

The main result of this paper is the following theorem, which provides an exact equation for the Hausdorff dimension of an average conformal hyperbolic set. 
\begin{thm} (Main theorem)
Let $f: M \rightarrow M$ be a hyperbolic diffeomorphism on a Riemannian manifold, with a locally maximal hyperbolic set $\Lambda$, and let $x \in \Lambda$. Suppose $\Lambda$ is average conformal. Then for any $x \in \Lambda$,
\begin{align*}
\dim_H(\Lambda \cap W^{(u)}(x))	&= \underline{\dim}_B (\Lambda \cap W^{(u)}(x)) = \overline{\dim}_B(\Lambda \cap W^{(u)}(x))\\
								&= \frac{h_{\kappa^{(u)}}(f) \dim E^{(u)}(x)}{\int_{\Lambda} \log |\det \left(d_x f|_{E^{(u)}} \right)| d\kappa^{(u)}},\\
\dim_H(\Lambda \cap W^{(s)}(x))	&= \underline{\dim}_B (\Lambda \cap W^{(s)}(x)) = \overline{\dim}_B(\Lambda \cap W^{(s)}(x))\\
								&= \frac{h_{\kappa^{(s)}}(f) \dim E^{(s)}(x)}{\int_{\Lambda} \log |\det \left(d_x f|_{E^{(s)}} \right)| d\kappa^{(s)}}.
\end{align*}
where $h_{\kappa^{(s)}}(f)$ and $h_{\kappa^{(s)}}(f)$ are the entropies of $f$ with respect to the unique equilibrium measures $\kappa^{(s)}$ and $\kappa^{(u)}$ corresponding to $\log |\det \left(d_x f|_{E^{(s)}} \right)|$ and $\log |\det \left(d_x f|_{E^{(u)}} \right)|$ respectively. 
\end{thm}
Furthermore, we prove as a corollary that $\dim_H \Lambda = \dim_H(\Lambda \cap W^{(s)}(x)) + \dim_H(\Lambda \cap W^{(u)}(x))$. 

\section{Preliminaries}

\subsection{Lyapunov exponents}
For a $C^1$ diffeomorphism $f: M \rightarrow M$ on a compact Riemannian manifold $M$, a point $x \in M$ and a nonzero vector $u \in T_x M$, the \textit{Lyapunov exponent} of $u$ is defined by
$$\chi(x, u) = \lim_{n \rightarrow \infty} \frac{1}{n} \log \|d_xf^n(u)\|,$$
if the limit exists. Vectors $u$ with the same Lyapunov exponent $\chi$ (plus the zero vector) form a linear subspace $E^\chi(x)$ of $T_x M$ called the Lyapunov space of $\chi$. These spaces form an invariant bundle in the sense that $T_x f^n(E^\chi(x)) = E^\chi(f^n x)$, for all $n \in \mathbb{Z}$. 

Let $\mathcal{E}$ be the set of all ergodic invariant measures on $M$. According to the Oseledets theorem \cite{Oseledets 1}, for any $x \in M$ there are finitely many Lyapunov exponents, $\chi_1(x) < \ldots < \chi_{l(x)}(x)$. Furthermore $l(x)$ and $\chi_i(x)$ are constant for $\mu$-almost every $x$, so for these values of $x$ we denote them by $l(\mu)$ and $\chi_i(\mu)$. This is not the same definition as the one used in \cite{average conformal 1}, but it is equivalent.

\subsection{Conformal, weakly conformal, quasi-conformal, and average conformal maps}

Let $M$ be a compact manifold and let $f: M \rightarrow M$ be a continuous (not necessarily differentiable) map. There are three generalizations of conformal repellers: weakly conformal, quasi-conformal and average conformal, any of which can be extended to hyperbolic sets by applying the same conditions to the stable and unstable manifolds seperately. Weakly conformal maps were first defined in \cite{Pesin}, quasi-conformal maps for induced expansive maps first appeared in \cite{Barreira 1} under the name ``asymptotically conformal maps'' and were extended to all continuous expansive maps and renamed quasi-conformal in \cite{Pesin}. 

When the map $f$ is $C^1$, average conformal is the most general condition; it is easy to show that a conformal repeller is weakly conformal, a weakly conformal repeller is quasi-conformal and a quasi-conformal repeller is average conformal. Furthermore, these definitions can be extended to hyperbolic sets, where these facts still hold. Average conformal maps also have the advantage of a simpler definition; it is relatively easy to check if the Lyapunov exponents are equal. Note that if $f$ is not differentiable, the concepts conformal and average conformal are not well-defined.

A $C^1$ map $f$ with a repeller is called \textit{conformal} if its derivative is a multiple of an isometry, i.e. there exists a continuous function $a(x)$ and an isometry $\mbox{Isom}_x$ such that $d_x f = a(x) \mbox{Isom}_x$ is a multiple of an isometry, that is $d_x f = a(x)$. 

A $C^1$ map $f$ with a hyperbolic set is called $u$-\textit{conformal} (respectively, $s$-\textit{conformal}) if there exists a continuous function $a^{(u)}(x)$ (respectively, $a^{(s)}(x)$), such that $d_x f = a^{(u)}(x) \mbox{Isom}_x$ (respectively, $d_x f = a^{(s)}(x) \mbox{Isom}_x$) for some isometry $\mbox{Isom}_x$. Then $f$ is called \textit{conformal} if it is both $u$-conformal and $s$-conformal.

We will not go into detail for weakly conformal maps here, except to say that they are defined for continuous maps (not necessarily $C^1$), and they are more general than conformal maps but less general than quasi-conformal maps. A definition can be found in \cite{Pesin} (page 191). 

Quasi-conformal maps (originally called asymptotically conformal in \cite{Barreira 1}, not to be confused with quasiconformal mappings) are defined on continuous expanding maps with a Markov partition $\{R_i\}_i$ (see e.g. \cite{Katok Hasselblatt} for a definition of Markov partitions). For a point $z \in \Lambda$, let $R(z)$ be the rectangle containing. Let $z \in \Lambda$, and let $k \geq 0$, $n \geq 1$. Then define numbers
$$\underline{\lambda}_k(z, n) = \inf_{C(x, n+k)} \left \{\frac{\|f^n x - f^n y\|}{\|x - y\|} \right\},$$
$$\overline{\lambda}_k(z, n) = \sup_{C(x, n+k)} \left \{\frac{\|f^n x - f^n y\|}{\|x - y\|} \right\},$$
where the infimum and the supremum are taken over a cylinder set $$\displaystyle C(z, n+k) = \bigcap_{j=0}^{n+k} f^{-j} R(f^j z).$$ This set becomes smaller as $n + k$ increases, i.e. $C(z, n + k + 1) \subset C(z, n + k)$. See \cite{Barreira 1} or \cite{Pesin} for more details.

\begin{defn} \cite{Pesin}
We say that the map is \textit{quasi-conformal} if there exist numbers $C > 0$ and $k > 0$ such that for all $x \in \Lambda$ and $n \geq 0$, 
$$\overline{\lambda}_k(\omega,n) \leq C \underline{\lambda}_k(\omega, n)$$
\end{defn}

Now we define quasi-conformal hyperbolic sets. Let $f: M \rightarrow M$ with a basic hyperbolic set $\Lambda$ with a Markov partition $\{R_i \}$. For a sequence $\omega = (\ldots, i_{-1}, i_0, i_1, \ldots)$ define numbers

$$\underline{\lambda}^{(s)}_k(z, n) = \inf_{C^{(s)}(z, n+k)} \left \{\frac{\|f^n x - f^n y\|}{\|x - y\|} \right\},$$
$$ \overline{\lambda}^{(s)}_k(z, n) = \sup_{C^{(s)}(z, n+k)} \left \{\frac{\|f^n x - f^n y\|}{\|x - y\|} \right\},$$
$$\underline{\lambda}^{(u)}_k(z, n) = \inf_{C^{(u)}(z, n+k)} \left \{\frac{\|f^n x - f^n y\|}{\|x - y\|} \right\},$$
$$ \overline{\lambda}^{(u)}_k(z, n) = \sup_{C^{(u)}(z, n+k)} \left \{\frac{\|f^n x - f^n y\|}{\|x - y\|} \right\},$$
where the infimums and supremums are taken over cylinder sets $$\displaystyle C^{(s)}(z, n + k) = \bigcap_{j=-(n+k)}^{n+k} f^{-j} R(f^j z) \cap W^{(s)}(z)$$ and $$\displaystyle C^{(u)}(z, n + k) = \bigcap_{j=-(n+k)}^{n+k} f^{-j} R(f^j z) \cap W^{(u)}(z).$$

\begin{defn}
We say that $f$ is $s$-\textit{quasi-conformal} if there exist numbers $C > 0$ and $k > 0$ such that for all $z \in \Lambda$ and $n \geq 0$, 
$$\overline{\lambda}^{(s)}_k(\omega,n) \leq C \underline{\lambda}_k(z, n).$$
We say that $f$ is $u$-\textit{quasi-conformal} if there exist numbers $C > 0$ and $k > 0$ such that for all $z \in \Lambda$ and $n \geq 0$, 
$$\overline{\lambda}^{(u)}_k(\omega,n) \leq C \underline{\lambda}_k(z, n).$$
Finally we say $f$ is \textit{quasi-conformal} if it is both $u$-quasi-conformal and $s$-quasi-conformal.
\end{defn}

Now we give a sketch of the proof that average conformal hyperbolic sets are more general than $C^1$ quasi-conformal hyperbolic sets.
\begin{thm}
If $f$ is a $C^1$ quasi-conformal hyperbolic map then $\Lambda$ is an average conformal hyperbolic set.

\begin{proof}
Let $f$ be a $u$-quasi-conformal hyperbolic map. Then there exists $k$ such that $\overline{\lambda}^{(u)}_k(z,n) \leq C \underline{\lambda}^{(u)}_k(x,n)$. For $k+1$, the infimum and supremum are taken over a smaller set, so we have for any $z$ and $n$,
$$\overline{\lambda}^{(u)}_{k+1}(z, n) \leq \overline{\lambda}^{(u)}_k(z, n) \leq C \underline{\lambda}^{(u)}_k(z,n) \leq C \underline{\lambda}^{(u)}_{k+1}(z,n).$$
For a matrix $A$, we define $m(A) = \|A^{-1}\|^{-1}$. If $f$ is $C^1$, then 
$$\lim_{k \rightarrow \infty} \inf_{C^{(u)}(z,n+k)} \left \{\frac{\|f^n x - f^n y\|}{\|x - y\|} \right\} = m \left( d_z f^n|_{E^{(u)}} \right),$$
$$\lim_{k \rightarrow \infty} \sup_{C^{(u)}(z,n+k)} \left \{\frac{\|f^n x - f^n y\|}{\|x - y\|} \right\} = \|d_z f^n|_{E^{(u)}} \|.$$
So for any regular point $z$ and any $u \in E^{(u)}(z)$, the Lyapunov exponents satisfy
\begin{align*}
\chi(z, u)	&= \lim_{n \rightarrow \infty} \frac{1}{n} \log \|\left(d_z f^n \right) u\|\\
			&\leq \lim_{n \rightarrow \infty} \lim_{k \rightarrow \infty} \frac{1}{n} \log \sup_{C^{(u)}(z,n+k)} \left \{\frac{\|f^n x - f^n y\|}{\|x - y\|} \right\}\\
			&\leq \lim_{n \rightarrow \infty} \lim_{k \rightarrow \infty} \frac{1}{n} \log \overline{\lambda}^{(u)}_k(z, n) \leq  \lim_{n \rightarrow \infty} \lim_{k \rightarrow \infty} \frac{1}{n} \log C \underline{\lambda}^{(u)}_k(z, n)\\
			&\leq  \lim_{n \rightarrow \infty} \frac{C}{n} + \frac{1}{n} \log m \left(d_z f^n|_{E^{(u)}} \right) = \inf_u \chi^{(u)}(z,u).
\end{align*}
Similarly, if $f$ is an $s$-quasi-conformal map, then for all $z \in \Lambda$ and $u \in E^{(s)}$ we have $\displaystyle \chi(z,u) \leq \inf_{u \in E^{(s)}} \chi(z,u)$. So if $f$ is quasi-conformal, there are exactly two Lyapunov exponents, $\chi^{(s)}$ and $\chi^{(u)}$, which implies $f$ is average conformal.
\end{proof}
\end{thm}

An example of an average conformal repeller that is not conformal can be found in \cite{example}. However that example is also quasi-conformal. It may be difficult to find an explicit example of an average conformal map that is not quasi-conformal (whether it has a repeller or a hyperbolic set). 

\subsection{Sub-additive and super-additive sequences}

From here on, we assume that $M$ is a compact manifold, and $f: M \rightarrow M$ is a $C^1$ map with a locally maximal hyperbolic set $\Lambda$. A sequence of continuous functions $\phi_n: X \rightarrow \mathbb{R}$ is called \textit{sub-additive} if 
$$\phi_{m+n}(x) \leq \phi_n(x) + \phi_m(f^n x),$$
and \textit{super-additive} if
$$\phi_{m+n}(x) \geq \phi_n(x) + \phi_m(f^n x).$$
We define the following four function sequences:
\begin{itemize}
	\item $\mathcal{F}^{(+, u)} = \{\Phi^{(u)}_n(x)\} = \{- \log \|d_x f^n|_{E^{(u)}}\| \}$ is a super-additive sequence.
	\item $\mathcal{F}^{(-, u)} = \{\varphi^{(u)}_n(x)\} = \{- \log m \left(d_x f^n|_{E^{(u)}} \right) \}$ is a sub-additive sequence.
	\item $\mathcal{F}^{(+, s)} = \{\Phi^{(s)}_n(x)\} = \{- \log \|d_x f^n|_{E^{(s)}}\| \}$ is a super-additive sequence.
	\item $\mathcal{F}^{(-, s)} = \{\varphi^{(s)}_n(x)\} = \{- \log m \left(d_x f^n|_{E^{(s)}} \right) \}$ is a sub-additive sequence.
\end{itemize}
For most of this paper it makes no difference whether we work on the stable or unstable manifold, so we will often set either $E = E^{(u)}$ or $E = E^{(s)}$, and then write $\mathcal{F}^+ = \{\Phi_n(x)\} = \{- \log \|d_x f^n|_E\| \}$ and $\mathcal{F}^- = \{\varphi_n(x)\} = \{- \log m \left( d_x f^n|_E \right) \}$.

\subsection{Topological Pressure}
A set $E \subset X$ is called $(f, n, \epsilon)$-separated with respect to $f$ if for every $x, y \in E$, $d_n(x, y) = \displaystyle \max_{0 \leq i \leq n-1} d(f^ix, f^i y) > \epsilon$. Define the \textit{Birkhoff sum} of a function $\phi$ with respect to $f$ by $S_n\phi(x) = S_n[f]\phi(x) = \sum_{j=0}^{n-1} \phi(f^j x)$. Define the topological pressure of a continuous function $\phi$ by
$$P_n(f, \phi, \epsilon) = \sup \left\{ \sum_{x \in E} \exp S_n\phi(x) : E \mbox{ is $(f,n,\epsilon)$-separated} \right\},$$
and
$$P(f, \phi) = \lim_{\epsilon \rightarrow 0} \limsup_{n \rightarrow \infty} \frac{1}{n} \log P_n(f, \phi, \epsilon).$$
Let $\mathcal{F} = \{\phi_n\}$ be a function sequence. Then the topological pressure of $\mathcal{F}$ is defined by
$$P^*_n(f, \mathcal{F}, \epsilon) = \sup \left\{ \sum_{x \in E} \exp \phi_n(x) : E \mbox{ is $(f,n,\epsilon)$-separated} \right\},$$
and
$$P^*(f, \mathcal{F}) = \lim_{\epsilon \rightarrow 0} \limsup_{n \rightarrow \infty} \frac{1}{n} \log P_n(f, \mathcal{F}, \epsilon).$$

Let $h_\mu(f)$ denote the entropy of $f$ with respect to a measure $\mu$. We use $P_B(f, \mathcal{F})$ to denote Barreira's definition of topological pressure for subadditive continuous functions via open covers \cite{Barreira 1}. 

\begin{prop} \cite{Thermodynamic formalism}
Assume that $h_\mu(f) < \infty$ and that the map $\mu \mapsto h_\mu(f)$ is upper-semi continuous. Then $P^*(f, \mathcal{F}) = P_B(f, \mathcal{F})$. 
\end{prop}

Let $\mathcal{M}(X)$ be the space of all Borel probability measures endowed with the weak* topology. Let $\mathcal{M}(X, f)$ be the subspace of $\mathcal{M}(X)$ consisting of $f$-invariant measures. Let $\mathcal{E}(f)$ be the set of ergodic $f$ invariant measures.

\section{Theorems for average conformal hyperbolic sets}

The following proposition is proven and used frequently in \cite{average conformal 1}, but not stated as a seperate result. 
\begin{prop} \label{4C prop}
Let $\{\phi_n(x)\}$ be a sub-additive sequence of functions on $M$, and fix some $m \in \mathbb{N}$. Then for any $n \in \mathbb{N}$, 

$$\phi_n(x) \leq \sum_{j=0}^{n-1} \frac{1}{m} \phi_m(f^j x) + 4C_1(m) = S_n\left(\frac{\phi_m}{m} \right)(x) + 4C_1,$$
where $C_1(m) = \displaystyle \max_{i = 1, \ldots, 2m-1} \max_{x \in M} \phi_i(x)$. Similarly, if $\{\phi_n(x)\}$ is a super-additive function sequence, then 
$$\phi_n(x) \geq \sum_{j=0}^{n-1} \frac{1}{m} \phi_m(f^j x) + 4C_2(m) = S_n\left(\frac{\phi_m}{m}\right)(x) + 4C_2,$$
where $C_2(m) = \displaystyle \min_{i = 1, \ldots, 2m-1} \min_{x \in M} \phi_i(x)$. 
\end{prop}

\begin{prop} \label{Kingman} (Kingman's sub-additive ergodic theorem)
Let $\phi_n$ be a sub-additive sequence of functions. Then
$$\lim_{n \rightarrow \infty} \frac{\phi_n}{n} = \inf_{n \geq 1} \frac{\phi_n}{n}.$$
\end{prop}

The following theorem is essentially the same as Theorem 4.2 in \cite{average conformal 1} for repellers. 
\begin{thm}
If $f$ is a diffeomorphism with a compact, hyperbolic, average conformal invariant set $\Lambda$, then for all $x \in \Lambda$ and any invariant bundle $E$ of $T_x M$, 

\begin{equation}\label{Fn = 0}
\displaystyle \lim_{n \rightarrow \infty} \frac{1}{n} \left(\log \|d_x f^n|_E\| - \log m \left(d_x f^n|_E \right) \right) = 0
\end{equation}
uniformly on $\Lambda$. 

\begin{proof}
The argument here is almost identical to that in \cite{average conformal 1} (Theorem 4.2). For $E = E^{(u)}$ or $E = E^{(s)}$, let $$F_n(x) = \log\|d_x f^n|_E\| - \log m(d_x f^n|_E), n \in \mathbb{N}, x \in \Lambda.$$ Suppose equation (\ref{Fn = 0}) is false. Then there exist sequences $n_k \geq k$ and $x_{n_k} \in \Lambda$ such that for all $k \geq 0$,
$$\frac{1}{n_k} F_{n_k}(x_{n_k}) \geq \epsilon_0.$$
Define measures
$$\mu_{n_k} = \frac{1}{n_k} \sum_{i=0}^{n_k -1} \delta_{f^i(x_{n_k})}.$$
Since $\mathcal{E}(f)$ is compact, there exists a subsequence of $\mu_{n_k}$ that converges to $\mu$. Without loss of generality suppose that $\mu_{n_k} \rightarrow \mu$. It is easy to show that $\mu$ is $f$-invariant. 
Then using Proposition \ref{4C prop} and the argument in \cite{average conformal 1}, we have
$$\lim_{m \rightarrow \infty} \int_M \frac{1}{m} F_m(x) d\mu \geq \epsilon_0 >0.$$
By the ergodic decomposition theorem (see Remark 2 in \cite{Walters}) there is an ergodic measure $\tilde{\mu}$ satisfying the same inequality.
By Kingman's subadditive ergodic theorem (Theorem 10.1 in \cite{Walters}), 
$$\lim_{m \rightarrow \infty} \frac{1}{m} \int_M \Phi_m d\tilde{\mu}= \int_M \chi d\mu = \chi \mbox{ and } \lim_{m \rightarrow \infty} \frac{1}{m} \int_M \varphi_m d\tilde{\mu} = \int_M \chi d\tilde{\mu} = \chi,$$
where $\chi$ is the Lyapunov exponent associated with $E$ (e.g. $\chi^{(u)}$ if $E = E^{(u)}$). So we have 
$$\lim_{m \rightarrow \infty} \int_M \frac{1}{m} F_m(x) \tilde{\mu} = 0,$$
which proves equation (\ref{Fn = 0}) by contradiction.
\end{proof}
\end{thm}

\subsection{Variational principle}
The following theorem unifies the variational principles for sub-additive sequences (Theorem 1.1 in \cite{Thermodynamic formalism}) and super-additive sequences (Theorem 5.1 in \cite{average conformal 1}).

\begin{thm} \label{variational principle} (Variational principle for subadditive and super-additive functions)
Let $f: X \rightarrow X$ be a hyperbolic diffeomorphism, let $E$ be an invariant subbundle of $T_xM$ and let $\mathcal{F}^+ = \{\Phi_n(x)\} = \{- \log \|d_x f^n|_E\| \}$, $\mathcal{F}^- = \{\varphi_n(x)\} = \{- \log m(d_x f^n|_E ) \}$. Then 
\begin{align}
P^*(f, \mathcal{F}^-)	&= \sup\{ h_\mu(f) + \lim_{n \rightarrow \infty} \frac{1}{n} \int_M \varphi_n	 d \mu : \mu \in \mathcal{M}(X, f) \}\label{VP 1}\\		
						&= \sup\{ h_\mu(f) + \lim_{n \rightarrow \infty} \frac{1}{n} \int_M \Phi_n d \mu : \mu \in \mathcal{M}(X, f) \} \label{VP 2}\\	
						&=P^*(f, \mathcal{F}^+). \label{VP 3}
\end{align}

\begin{proof}
Equation (\ref{VP 1}) is proven in \cite{Thermodynamic formalism}. Theorem \ref{Fn = 0} implies that
$$\lim_{n \rightarrow \infty} \frac{1}{n} \int_M \varphi_n d \mu = \lim_{n \rightarrow \infty} \frac{1}{n} \int_M \Phi_n d \mu,$$
which gives equation (\ref{VP 2}). For a fixed $m$, let $n = mk + l$, $0 \leq l < m$. Since $\{\Phi_n\}$ is super-additive, Proposition \ref{4C prop} gives
$$\Phi_n(x) \geq \sum_{j=0}^{n-1} \frac{1}{m} \Phi_m(f^j x) + 4C_2.$$
Following the argument in Theorem 5.1 of \cite{average conformal 1}, we have 

\begin{align*}
P^*(f, \mathcal{F}^+)	&\geq \sup \left\{ h_\mu(f) + \lim_{n \rightarrow \infty} \frac{1}{n} \int_M \Phi_n d \mu : \mu \in \mathcal{M}(X, f) \right\}\\
						&=P^*(f, \mathcal{F}^-)
\end{align*}
But $\Phi_n(x) \leq \phi_n(x)$ implies $P^*(f, \mathcal{F}^+) \leq P^*(f, \mathcal{F}^-)$, which gives equation (\ref{VP 3}).

\end{proof}
\end{thm}

The following lemma corresponds to Lemma 6.1 in \cite{average conformal 1}. The proof is similar, but we go into more detail for the pressure. 
\begin{lem} \label{lemma A}
If $\phi_n(x)$ is a subadditive sequence, then
$$\lim_{k \rightarrow \infty} \frac{1}{2^k} P(f^{2^k}, \phi_{2^k}) \leq \lim_{m \rightarrow \infty} P(f, \frac{\phi_{2^m}}{2^m}).$$

\begin{proof}

For a fixed $m < k$, let $C_1 = C_1(m) = \displaystyle \max_{x \in M} \max_{j = 1, \ldots 2^m} \phi_j(x)$. Then by Proposition \ref{4C prop}, 
$$\phi_{2^k}(f^{2^k l}x) \leq \displaystyle \sum_{j = 0}^{2^k - 1} \frac{1}{2^m} \phi_{2^m}(f^j f^{2^k l} x) + 4C_1.$$
The Birkhoff sums satisfy
$$S_n[f^{2^k}] \phi_{2^k}(x) \leq S_{n 2^k}[f]\left(\frac{\phi_{2^m}}{2^m}\right)(x) + 4nC_1$$
This means that 
\begin{align*}
P_n(f^{2^k}, \phi_{2^k}, \epsilon)	&=   \sup_E \left\{ \sum_{x \in E} \exp S_n[f^{2^k}] \phi_{2^k}(x) \right\}\\
								 	&\leq\sup_E \left\{ \sum_{x \in E} \exp \left( S_{n 2^k}[f]\left(\frac{1}{2^m}\phi_{2^m} \right)(x)+ 4nC_1 \right) \right\},
\end{align*}
where the supremum is over $(f^{2^k}, n, \epsilon)$-separated subsets $E$. For a fixed $k \in \mathbb{N}$, if $E \subset M$ is an $(f^{2^k}, n, \epsilon)$-separated set, then $E$ is also an $(f, n 2^k, \epsilon)$-separated set. So
\begin{align*}
P(f^{2^k}, \phi_{2^k})	&\leq \lim_{\epsilon \rightarrow 0} \limsup_{n \rightarrow \infty} \frac{1}{n} \log \left( P_{n 2^k}\left(f, \frac{1}{2^m}\phi_{2^m}, \epsilon \right) e^{4nC_1} \right)\\
						&\leq 2^k P\left(f, \frac{1}{2^m}\phi_{2^m} \right) + 4C_1.
\end{align*}
Thus, for all $\displaystyle m \in \mathbb{Z}^+$, $\displaystyle \lim_{k \rightarrow \infty} \frac{1}{2^k} P(f^{2^k}, \phi_{2^k})	\leq P(f, \frac{1}{2^m} \phi_{2^m})$, so in particular,
$$\lim_{k \rightarrow \infty} \frac{1}{2^k} P(f^{2^k}, \phi_{2^k}) \leq \lim_{m \rightarrow \infty} P\left(f, \frac{\phi_{2^m}}{2^m} \right)$$ 
and Proposition \ref{Kingman} implies that the limit on the right exists.
\end{proof}
\end{lem}

The map $f$ is \textit{expansive} if there exists $\varepsilon > 0$ such that for any $x \neq y \in \Lambda$, there exists $n \in \mathbb{Z}$ such that $d(f^nx, f^n y) \geq \varepsilon$. Any repeller or hyperbolic diffeomorphism is expansive. It is well known (see e.g. Theorem 8.2 in \cite{Walters}) that if $f$ is an expansive homeomorphism, the entropy map $h_\mu(f)$ is upper semi-continuous with respect to $\mu$. 

The next lemma corresponds to Lemma 6.2 in \cite{average conformal 1}. 
\begin{lem} \label{lemma B}
If $\phi_n(x)$ is a sub-additive sequence, then
$$\lim_{k \rightarrow \infty} P\left(f, \frac{\phi_{2^k}}{2^k} \right) \leq P^*(f, \{\phi_n\}).$$

\begin{proof}
By the variational principle, for any $k \in \mathbb{Z}^+$ there exists $\mu_{2^k} \in \mathcal{M}(f|_\Lambda)$ such that
$$P\left(f, \frac{\phi_{2^k}}{2^k} \right) = h_{\mu_{2^k}}(f) + \int_{\Lambda} \frac{\phi_{2^k}}{2^k} d\mu_{2^k}.$$
Since $\mathcal{M}(f|_\Lambda)$ is compact, $\mu_{2^k}$ has a subsequence that converges to $\mu \in \mathcal{M}(f|_\Lambda)$. Without loss of generality, suppose that $\mu_{2^k}$ converges to $\mu$. Fix some $s \in \mathbb{N}$. Then from the sub-additivity of $\phi$ and the invariance of $\mu_{2^k}$, we have for all $k > s$, 

\begin{align*}
\int_{\Lambda} \frac{\phi_{2^k}(x)}{2^k}  d\mu_{2^k}	&\leq \int_{\Lambda} \frac{\phi_{2^k}(x)}{2^k}  d\mu_{2^k}\\
													&\leq \int_{\Lambda} \frac{1}{2^k} (\phi_{2^s}(x) + \phi_{2^s}(f^{2^s}x) + \ldots + \phi_{2^s}(f^{2^k}x)) d\mu_{2^k}\\
													&\leq \int_{\Lambda} \frac{2^{k-s}}{2^k} \phi_{2^s}(x) d\mu_{2^k} = \int_{\Lambda} \frac{\phi_{2^s}(x)}{2^s} d\mu_{2^k}.
\end{align*}
Since $h_\mu(f)$ is upper semi-continuous with respect to $\mu$, we have 
\begin{align*}
\lim_{k \rightarrow \infty} P(f, \frac{\phi_{2^k}}{2^k})	&= \lim_{k \rightarrow \infty} \left( h_{\mu_{2^k}}(f) + \int_{\Lambda} \frac{\phi_{2^k}(x)}{2^k} d\mu_{2^k}\right)\\
															&\leq h_\mu(f) + \int_{\Lambda} \frac{\phi_{2^s}(x)}{2^s} d\mu
\end{align*}
for all $s \in \mathbb{N}$. So by the variational principle \ref{variational principle}, 
$$\lim_{k \rightarrow \infty} P(f, \frac{\phi_{2^k}}{2^k}) \leq h_\mu(f) + \lim_{s \rightarrow \infty} \int_{\Lambda} \frac{\phi_{2^s}(x)}{2^s} d\mu \leq P^*(f, \mathcal{F}).$$
Proposition \ref{Kingman} implies that the limit on the right exists. 

\end{proof}
\end{lem}

\begin{prop} (First part of the proof of Theorem 6.2 in \cite{average conformal 1})
\label{Ps increasing}
For all $s \geq 0$, the sequence $\frac{1}{2^k} P(f^{2^k}, s \Phi_{2^k})$ is monotone increasing in $k$.
\begin{proof}

The Birkhoff sum $S_n \Phi_{2^{k+1}}$ with respect to $f^{2^{k+1}}$ has the following property:
$$S_n[f^{2^{k+1}}] \Phi_{2^{k+1}}(x) = S_{2n}[f^{2^k}] \Phi_{2^k}(x).$$

Since $f$ is uniformly continuous, for all $\epsilon > 0$ there exists $\delta > 0$ such that if $E \subset M$ is an $(f^{2^{k+1}},n, \epsilon)$-separated set then $E$ is a $(f^{2^k},2n, \delta)$-separated set, and $\delta \rightarrow 0$ when $\epsilon \rightarrow 0$. So the pressure satisfies
\begin{align*}
P_n(f^{2^{k+1}},  s \Phi_{2^{k+1}}, \epsilon)	&=	\sup \left\{\begin{array}{lr}
       \displaystyle \sum_{x \in E} \exp(S_n[f^{2^{k+1}}] s\Phi_{2^{k+1}}):\\
       \mbox{$E$ is $(f^{2^{k+1}}, n, \epsilon)$-separated}
     \end{array}  \right \}\\
												&\geq	\sup_E \left\{\begin{array}{lr}
												\displaystyle \sum_{x \in E} \exp(S_{2n}[f^{2^k}] s\Phi_{2^k}):\\
												\mbox{$E$ is $(f^{2^k}, 2n, \delta)$-separated} \end{array} \right\}\\
												&=		P_{2n}(f^{2^k}, s\Phi_{2^k}, \delta),
\end{align*}
and
\begin{align*}
P(f^{2^{k+1}}, s\Phi_{2^{k+1}})	&= \lim_{\epsilon \rightarrow 0} \limsup_{n \rightarrow \infty} \frac{1}{n} P_n(f^{2^{k+1}}, s\Phi_{2^{k+1}}, \epsilon)\\
									&\geq 2\lim_{\delta \rightarrow 0} \limsup_{n \rightarrow \infty} \frac{1}{2n} P_{2n}(f^{2^k}, s\Phi_{2^k}, \delta)\\
									&\geq 2P(f^{2^k}, s\Phi_{2^k}).
\end{align*}
Therefore $\frac{1}{2^k} P(f^{2^k}, s \Phi_{2^k})$ is monotone increasing.

\end{proof}
\end{prop}

\begin{prop} (First part of the proof of Theorem 6.3 in \cite{average conformal 1})
\label{Pt decreasing}
For all $t \geq 0$, the sequence $\frac{1}{2^k} P(f^{2^k}, t \varphi_{2^k})$ is monotone decreasing in $k$.
\begin{proof}
The argument is very similar to the previous proposition.
\end{proof}
\end{prop}

\begin{prop} (Second part of the proof of Theorem 6.2 in \cite{average conformal 1})
\label{1/k s}
For all $k \in \mathbb{N}$, 
$$P^*(f, \mathcal{F}^+) \geq \frac{1}{k} P(f^k, \Phi_k).$$

\begin{proof}
For a fixed $k \in \mathbb{N}$, let $n = km + r$ with $0 \leq r < k$, and $C(k) = \displaystyle \min_{x \in M} \min_{1 \leq j \leq k} \Phi_j(x)$. Since $f$ is uniformly continuous, for $\epsilon > 0$ there exists $\delta > 0$ such that if $E \subset M$ is $(f, n, \epsilon)$-separated then it is $(f^k, m, \delta)$-separated, and $\delta(\epsilon) \rightarrow 0$ as $\epsilon \rightarrow 0$. Using the super-additivity of $\Phi_n$ we have
$$\Phi_n(x) \geq \Phi_k(x) + \Phi_k(f^k x) + \ldots + \Phi_k(f^{(m-1)k}) x) + \Phi_r(f^{mk}x).$$
Thus,
\begin{align*}
P_n^*(f, \mathcal{F}, \epsilon)	&=	  \sup \left\{ \sum_{x \in E} e^{\Phi_n(x)} : E \mbox{ is $(f, n, \epsilon)$-separated}\right\}\\
								&\geq \sup \left\{ \sum_{x \in E} e^{S_m \Phi_k(x)} e^{\Phi_r(f^{mk}x)} : E \mbox{ is $(f^k, m,\delta)$-separated} \right\}\\
								&\geq P_m(f^k, \Phi_k, \delta) e^C.
\end{align*}
So,
$$\frac{1}{n} \log P_n^*(f, \mathcal{F}^+, \epsilon) \geq \frac{1}{k m + r} \log P_m(f^k, \Phi_k, \delta) + \frac{1}{n}C.$$
Taking limits, we have
$$\lim_{\epsilon \rightarrow \infty} \limsup_{n \rightarrow \infty} \frac{1}{n} \log P_n^*(f, \mathcal{F}^+, \epsilon) \geq \frac{1}{k} \lim_{\delta \rightarrow \infty} \limsup_{m \rightarrow \infty} \frac{1}{m} \log P_m(f^k, \Phi_k, \delta),$$
which by definition gives
$$P^*(f, \mathcal{F}^+) \geq \frac{1}{k}P(f^k, \Phi_k).$$ 
\end{proof}
\end{prop}

\begin{prop} (Second part of the proof of Theorem 6.3 in \cite{average conformal 1})
\label{1/k t}
For all $k \in \mathbb{N}$, 
$$P^*_n(f, \mathcal{F}^-) \leq \frac{1}{k} P(f^k, \varphi_k).$$

\begin{proof}
The argument is very similar to the previous proposition.
\end{proof}
\end{prop}

\section{Proof of main theorem}
We will use the following well known theorem of Barreira (Theorem 3.18 in \cite{Barreira 1}).
\begin{thm} \label{Barreira thm} \cite{Barreira 1}
Let $f$ be a hyperbolic diffeomorphism. Let $s_n$ and $t_n$ be the unique roots of the Bowen equations $P(f^n, -t \log\|d_x f^n|_{E^{(u)}}\|) = 0$ and $P(f^n, -t \log m\left( d_x f^n|_{E^{(u)}} \right)) = 0$ respectively. Then for any $x \in \Lambda$,
$$s_n \leq \dim_H(\Lambda \cap W^{(u)}(x)) \leq \underline{\dim}_B (\Lambda \cap W^{(u)}(x)) \leq \overline{\dim}_B(\Lambda \cap W^{(u)}(x)) \leq t_n.$$

Similarly, if $s_n$ and $t_n$ are the unique roots of $P(f^n, -t \log\|d_x f^n|_{E^{(s)}}\|) = 0$ and $P(f^n, -t \log m\left( d_x f^n|_{E^{(s)}} \right)) = 0$, then for any $x \in \Lambda$,
$$s_n \leq \dim_H(\Lambda \cap W^{(s)}(x)) \leq \underline{\dim}_B (\Lambda \cap W^{(s)}(x)) \leq \overline{\dim}_B(\Lambda \cap W^{(s)}(x)) \leq t_n.$$
\end{thm}

Since the proofs for the stable and unstable components are identical, from now on we will use $E$ to denote either $E^{(s)}$ or $E^{(u)}$. 

\begin{thm} (Theorem 6.2 in \cite{average conformal 1}) \label{s theorem}
The sequence $\{s_{2^k}\}$ is monotone increasing and $s_{2^k} \rightarrow s_*$ as $k \rightarrow \infty$, where $s_*$ is the root of the equation $$P^*(f, -s_* \log \|d_x f^n|_E\|) = 0.$$

\begin{proof}
First we show that $\{s_{2^k}\}$ is monotone increasing in $k$. The function $P(f^{2^k}, s \phi_{2^k})$ is monotone decreasing in $s$ and by Proposition \ref{Ps increasing} it is monotone increasing in $k$, so its zero $s_{2^k}$ is monotone increasing in $k$. Hence the limit $\displaystyle \lim_{k \rightarrow \infty} s_{2^k}$ exists and we denote it by $\overline{s}$.

By Proposition \ref{Ps increasing}, we have $P(f^{2^{k+1}}, \Phi_{2^{k+1}}) \geq 2P(f^{2^k}, \Phi_{2^k})$ for all $k$. So if $s_{2^{k+1}}$ is the unique root of $P(t \Phi_{2^{k+1}})=0$, then
$$0 = P(f^{2^{k+1}}, s_{2^{k+1}} \Phi_{2^{k+1}}) \geq 2P(f^{2^k},  s_{2^{k+1}}\Phi_{2^k}).$$
By Proposition \ref{1/k s}, for all $k \in \mathbb{N}$ we have
$$P^*(f, s_{2^k} \mathcal{F}) \geq \frac{1}{2^k}P(f^{2^k}, s_{2^k}\Phi_{2^k}) = 0.$$
Next we show $P^*(f, \overline{s} \mathcal{F}) \leq 0$. For a fixed $m$, $$\frac{1}{2^m} P(f^{2^m}, s_{2^m} \Phi_{2^m}) = 0.$$
In the following we use the variational principle twice and the fact that $h_\mu(f) = \frac{1}{2^m} h_\mu(f^{2^m})$. There exists $\mu \in \mathcal{M}(f) \subset \mathcal{M}(f^{2^m})$ such that
\begin{align*}
P^*(f, \overline{s} \mathcal{F})	&= \left(h_\mu(f) + \overline{s} \lim_{m \rightarrow \infty} \frac{1}{2^m}  \int_M \Phi_{2^m} d\mu \right)\\
									&= \lim_{m \rightarrow \infty} \frac{1}{2^m} \left(h_\mu(f^{2^m}) 		  + s_{2^m} \int_M \Phi_{2^m} d\mu \right)\\
									&\leq \lim_{m \rightarrow \infty} \frac{1}{2^m} P(f^{2^m}, s_{2^m} \Phi_{2^m}) = 0.
\end{align*}
Since $P^*(f, \overline{s} \mathcal{F}) \geq 0$ and $P^*(f, \overline{s} \mathcal{F}) \leq 0$, we have $\overline{s} = s_*$ as required.
\end{proof}
\end{thm}

\begin{rem}
In \cite{average conformal 1} the authors incorrectly state that $s_{2^k}$ is monotone decreasing, however their proof implies that it is monotone increasing and they later assume it to be monotone increasing.
\end{rem}

\begin{thm} (Theorem 6.2 in \cite{average conformal 1}) \label{t theorem}
The sequence $\{t_{2^n}\}$ is monotone decreasing and 
$$\lim_{n \rightarrow \infty} t_{2^n} = t^*,$$
where $t^*$ is the unique root of $P(f, -t \log m(d_x f^n|_E))=0$.

\begin{proof}
First we show that $\{t_{2^n}\}$ is monotone decreasing. By Proposition \ref{Pt decreasing}, the function $P(f^{2^k}, t \phi_{2^k})$ is monotone decreasing in $t$ and monotone decreasing in $k$, so its zero $t_{2^k}$ is monotone decreasing in $k$. Hence the limit $\displaystyle \lim_{k \rightarrow \infty} t_{2^k}$ exists and we denote it by $\overline{t}$. By Proposition \ref{1/k t}, $$P(f, t \mathcal{F}^-) \leq \frac{1}{2^k} P(f^{2^k}, t \varphi_{2^k})$$ for all $k \in \mathbb{N}$. But by Lemmas \ref{lemma A} and \ref{lemma B}, we have
$$\lim_{k \rightarrow \infty} \frac{1}{2^k} P(f, t \varphi_{2^k}) \leq P^*(f, t\mathcal{F}^-).$$
So for all $t \geq 0$, 
$$\lim_{k \rightarrow \infty} \frac{1}{2^k} P(f, t \varphi_{2^k}) = P^*(f, t\mathcal{F}^-).$$
In particular, 
$$P^*(f, \overline{t} \mathcal{F}^-) = \lim_{k \rightarrow \infty} \frac{1}{2^k} P(f, t_{2^k} \varphi_{2^k}) = 0.$$
So $\overline{t} = t^*$ as required.

\end{proof}
\end{thm}

\begin{thm} \label{final theorem}
Let $\Lambda$ be a locally maximal hyperbolic set for a $C^1$ diffeomorphism $f$. Let $d^{(s)}$, $d^{(u)}$ be the dimensions of $E^{(s)}$ and $E^{(u)}$ respectively. Then
\begin{align*}
\dim_H(\Lambda \cap W^{(u)})	&= \underline{\dim}_B (\Lambda \cap W^{(u)}) = \overline{\dim}_B(\Lambda \cap W^{(u)}) = r^{(u)},\\
\dim_H(\Lambda \cap W^{(s)})	&= \underline{\dim}_B (\Lambda \cap W^{(s)}) = \overline{\dim}_B(\Lambda \cap W^{(s)}) = r^{(s)},
\end{align*}
where $r^{(u)} = \frac{h_{\kappa^{(u)}}(f)}{\int_{\Lambda} \frac{1}{d^{(u)}} \log |\det \left(d_x f|_{E^{(u)}} \right)| d\kappa^{(u)}}$ and $r^{(s)} = \frac{h_{\kappa^{(s)}}(f)}{\int_{\Lambda} \frac{1}{d^{(s)}} \log |\det \left(d_x f|_{E^{(s)}} \right)| d\kappa^{(s)}}$. 

\begin{proof}
Note that if $E$ is a subspace of $T_xM$ with dimension $d$, then 
$$\log m(d_x f^n|_E) \leq \frac{1}{d} \log |\det \left(d_x f^n|_E \right)| \leq \log \|d_x f^n|_E\|.$$
Also note that the sequence $\{\phi_n\} = \{ \frac{1}{d}\log |\det \left(d_x f^n|_E \right)|\}$ is additive, so $\phi_n = S_n \phi$ where $\phi = \frac{1}{d}\log |\det \left(d_x f|_E \right)|$. Therefore, by the definitions of pressure and the variational principle, there exists a measure $\kappa$ such that for any $t \in \mathbb{R}$, 
\begin{align*}
P^*(-t \frac{1}{d} \{\log |\det \left(d_x f^n|_E \right)|\}) &= P(-t \log |\det \left(d_x f|_E \right)|)\\
															 &= h_\kappa(f) - t \int_{\Lambda} \frac{1}{d}\log |\det \left(d_x f^n|_E \right)| d\kappa.
\end{align*}
So the solution $r$ of Bowen's equation $P^*(-r^{(u)} \{\frac{1}{d} \log |\det \left(d_x f^n|_E \right)|\}) = 0$ is given by 
$$r = \frac{h_{\kappa}(f)}{\int_{\Lambda} \frac{1}{d} \log |\det \left(d_x f|_E \right)| d\kappa}.$$ Setting $E = E^{(u)}$ or $E = E^{(s)}$ gives the result. 
\end{proof}

\end{thm}

\subsection{Dimension product structure}

For a point $x$ with an open neighbourhood $R$, let $h: W^{(u)}(x) \times W^{(s)}(x) \rightarrow R$ be the local product map. This map is always H\"{o}lder continuous, that is $d(h(x), h(y)) \leq C d(x,y)^\alpha$ for some $\alpha > 0$ and $C > 0$. If the diffeomorphism $f$ is conformal, the holonomy is Lipschitz. In \cite{Hasselblatt}, a point $x \in \Lambda$ is called $\alpha$-bunched if $\chi^{(s)}_{\min} - \chi^{(u)}_{\min} \leq \min\{\chi^{(s)}_{\max}, \chi^{(u)}_{\max}\}^{-\alpha}$. If $\Lambda$ is average conformal, then $\chi^{(s)}_{\max} = \chi^{(s)}_{\min}$ and $\chi^{(u)}_{\max} = \chi^{(u)}_{\min}$, so every point is at least $1$-bunched, and therefore Lipschitz continuous. It is well known that when the holonomy is Lipschitz, the dimension of the hyperbolic set is the sum of the dimensions of its stable and unstable manifolds. So we have the following corollary:

\begin{cor}
If $\Lambda$ is an average conformal hyperbolic set, then 
$$\dim_H \Lambda = \underline{\dim}_B \Lambda = \overline{\dim}_B \Lambda = r^{(s)} + r^{(u)},$$
where $r^{(s)}, r^{(u)}$ are defined as in Theorem \ref{final theorem}. 

\end{cor}

\end{document}